\documentclass[twocolumn]{autart}

\usepackage{graphicx}
\usepackage{cite}
\usepackage{comment}
\usepackage{amsmath,amssymb,amsfonts}
\usepackage{dcolumn}
\usepackage{bm}
\usepackage[normalem]{ulem}
\usepackage{comment}
\usepackage{float}

\usepackage{diagbox}

\usepackage{mathtools}
\usepackage{multirow}
\usepackage[utf8]{inputenc}
\usepackage{gensymb}
\usepackage{algorithm}
\usepackage{algpseudocode}
\usepackage{xcolor}

\usepackage[shortlabels]{enumitem}

\theoremstyle{plain}

\newtheorem{proposition}{Proposition}

\theoremstyle{definition}
\newtheorem{definition}{Definition}
\newtheorem{remark}{Remark}

\graphicspath{{./Figs/}}

\newcolumntype{C}[1]{>{\centering\arraybackslash$}m{#1}<{$}}
\newlength{\mycolwd}                                         
\newcommand{\ZA}[1]{{\color{black}{#1}}}

\def \bF {\pmb{F}}
\def \bf {\pmb{f}}
\def \bg {\pmb{g}}
\def \bh{\pmb{h}}
\def \bH {\pmb{H}}
\def \bD {\pmb{D}}

\def \bg {\pmb{g}}

\def \bu {\pmb{u}}

\def \bx {\pmb{x}}
\def \by {\pmb{y}}
\def \bX {\pmb{X}}

\def \bz {\pmb{z}}

\def \bs {\pmb{s}}

\def \vec {\textnormal{vec}}

\def \vec {\textnormal{vec}}

\def \barR {\bar{R}_{\|\cdot\|}}
\def \limR {R^{(\infty)}_{\|\cdot\|}}
\def\meanR {R_{\| \cdot \|}^{(N)}} 

\def \mean {\text{mean}}

\edef\endfrontmatter{%
  \unexpanded\expandafter{\endfrontmatter}
  \noexpand\endNoHyper 
}

\begin{document}

\begin{frontmatter}

\title{Bridging the Gap between Reactivity, Contraction,\\
 and Finite-Time Lyapunov Exponents}

\thanks[footnoteinfo]{This work was funded by NSF award 2037828 and The Simons Foundation MPS--TSM--00008005 (ZA).} 

\author[nm]{Amirhossein Nazerian}\ead{anazerian@unm.edu},    
\author[nm]{Francesco Sorrentino}\ead{fsorrent@unm.edu}$^1$\thanks[footnoteinfo]{Corresponding author.},  
\author[ia]{Zahra Aminzare}\ead{zahra-aminzare@uiowa.edu}      

\address[nm]{Mechanical Engineering Department, University of New Mexico, Albuquerque, NM, 87131 USA}                                         
\address[ia]{The Department of Mathematics, University of Iowa, Iowa City, IA 52242 USA}

\begin{keyword}                         
Reactivity; contraction theory; finite-time Lyapunov exponent.           
\end{keyword}                             
                      
\begin{abstract}
Reactivity, contractivity, and Lyapunov exponents are powerful tools for studying the stability properties of dynamical systems and have been extensively investigated in the literature for decades. In this paper, we review and extend the concepts of reactivity, contractivity, and finite-time Lyapunov exponents for discrete-time dynamical systems and establish connections among them. We focus on time-invariant maps, time-varying linear maps, and certain classes of time-varying nonlinear maps. In particular, we show that if the corresponding $p$-iteration systems (with $p > 1$) are contractive, then the original systems admit stable attractors such as fixed points or limit cycles. We demonstrate the application of these results to the analysis of synchronization stability in coupled networks and discuss how $p$-iteration systems can serve as a useful framework for studying network synchronization.
\end{abstract}

\end{frontmatter}

\section{Introduction}
The concept of reactivity was first introduced in \cite{Neubert1997ALTERNATIVES} for linear time-invariant systems. 
Reactivity measures the maximum amplification rate of the state norm at a given time \cite{Neubert1997ALTERNATIVES}.
Several studies expanded the concept of reactivity to networked systems with applications such as single-integrator consensus problem \cite{nazerian2023single}, linear discrete-time systems \cite{nazerian2023reactability}, 
and more \cite{Farrell1996Generalized,Hennequin2012Nonnormal,Tang2014Reactivity, Biancalani2017Giant,Asllani2018Structure,MUOLO2019Patterns,Gudowska2020From,Lindmark2021Centrality,Duan2022Network,nazerian2024efficiency}. 
All of the aforementioned references study the reactivity either associated with the linear or linearized dynamics based on the $L_2$ norm.
Here, we generalize the concept of reactivity to study the transient and steady-state behavior of discrete-time nonlinear systems using 
general norms.
Also, we show that this analysis, which is based on reactivity, closely relates to contraction theory and finite-time Lyapunov exponents.

Contraction theory studies the incremental convergence of the trajectories of dynamical systems \cite{lohmiller1998contraction,lohmiller2000nonlinear,slotine2001modularity,slotine2003modular,wang2005partial}, and can be used to assess  stability 
in terms of the time evolution of the distance between any two trajectories of a system.
If the system is contractive, the distance between the two trajectories diminishes in time, meaning that regardless of the initial condition, the trajectories will eventually converge to one another with a known finite convergence rate. Hence, contraction provides a more conservative convergence criterion than the conventional asymptotic stability.
Recent studies investigated contraction theory in various disciplines~\cite{bullo2022contraction,aminzare2013logarithmic,aminzare2014contraction}.

Lyapunov exponents are a widely used method for determining the asymptotic stability of linear time-varying systems across various fields, including synchronization \cite{nazerian2022matryoshka, Nazerian2023epl, nazerian2023commphys, panahi2022pinning, panahi2021group} and chaos detection \cite{Wolf1986Quantifying, abarbanel1991lyapunov}. 
{\color{black}
They measure the asymptotic rate of convergence or divergence of nearby trajectories.
On the other hand, the finite-time Lyapunov exponents (FTLEs) quantify the rate of growth of the norm of a perturbation over a finite time horizon and are particularly useful for analyzing the stretching and folding characteristics of attractors and other invariant sets \cite[Chapter 4.4]{ott2002chaos}.}

The stability analysis conducted for linearized systems using the Lyapunov exponents or reactivity generally guarantees local stability. 
For instance, the master stability function approach uses Lyapunov exponents to characterize the local stability of the synchronous solution in networks of coupled oscillators \cite{nazerian2022matryoshka, Nazerian2023epl, nazerian2023commphys, panahi2022pinning, panahi2021group,panahi2021cluster}.  
On the other hand, contraction theory 
yield conditions for global stability. For example, for diffusively coupled systems, contraction theory establishes conditions that ensure the existence and stability of a complete synchronous solution~\cite{aminzare2014contraction}.

In this work, we focus on time-invariant maps, time-varying linear maps, and certain classes of time-varying nonlinear maps. We review the notions of reactivity, contractivity, and Lyapunov exponents, and establish connections among them. In particular, for the purposes of this paper, the terms `reactive' and `non-contractive' as well as `contractive' and `non-reactive' are used interchangeably.
Our main contributions are threefold:
(1) We introduce the concept of \emph{sup-mean contractivity/reactivity} and show that, to ensure contractivity, it is sufficient for the contractivity condition to hold on average over a finite number of time steps, rather than at each individual step.
(2) We derive conditions for the existence and stability of attractors—such as fixed points and limit cycles—for discrete-time dynamical systems, based on the sup-mean contractivity of their $p$-iteration systems (with $p > 1$). We see that increasing the number of iterations $p$ enlarges the range of parameters for which stability is guaranteed.
(3) We apply our results to the study of synchronization in discrete-time dynamical systems over a synchronous attractor, and obtain conditions for the local stability of the synchronous solution based on the average transverse reactivity along this attractor.

The rest of the paper is organized as follows. Sec.~\ref{sec:pre} reviews the notations and preliminaries. 
Sec.~\ref{sec:p-th interation} introduces the concept of sup-mean contractivity/reactivity for discrete-time dynamical systems and establishes conditions for the existence and stability of fixed points and limit cycles, based on the sup-mean contractivity of their $p$-iteration systems.
Sec.~\ref{sec:cle} establishes a connection with the study of the Lyapunov exponents of a system.   Sect.~\ref{sec:sync} presents examples of application of the theory to the synchronization of networks of coupled discrete-time systems, and, finally, the conclusions are given in Sec.~\ref{sec:conc}.

\section{Preliminaries} \label{sec:pre}
In this section, we first define the reactivity of a map with respect to a general norm and then introduce the reactivity of discrete-time dynamical systems. The goal is to connect the concept of reactivity to the stability of the discrete dynamical system.  

\begin{definition} \normalfont{\textbf{(Reactivity of a map with respect to a general norm)}} \label{def:Reactivity1}
Let $\bf: \mathbb{R}^{m} \rightarrow \mathbb{R}^{m}$ be a Lipschitz continuous function and $\| \cdot \|$ be an arbitrary  norm in $\mathbb{R}^{m}$. 
The Lipschitz constant for function $\bf$ induced by $\| \cdot \|$ is defined as
$
    L_{\| \cdot \|} [\bf] := \sup_{\bx \neq \by} \frac{\| \bf(\bx) - \bf(\by) \|}{\| \bx - \by \|}.
$
The reactivity of $\bf$ induced by $\| \cdot \|$ is denoted by $r_{\| \cdot \|} [\bf]$ and is defined by
$ r_{\| \cdot \|} [\bf] := L_{\| \cdot \|} [\bf] - 1.$
    \end{definition}
Let $Q$ be an $m\times m$ positive definite matrix and $A$ be an arbitrary $m\times m$ matrix. Then, 
the reactivity induced by $Q-$weighted $L_1$, $L_2$, and $L_{\infty}$ norms for linear maps $\bf(\bx) = A \bx$ becomes:
\begin{align*}
&r_{1,Q} [A] = \max_{1\leq j \leq m} \sum_{i = 1}^m | (Q A Q^{-1})_{ij} | - 1; \\  
& r_{2,Q} [A] = \sigma_1(Q A Q^{-1}) - 1;\\  
& r_{\infty,Q} [A] = \max_{1\leq i \leq m} \sum_{j = 1}^m | (Q A Q^{-1})_{ij} | - 1.
\end{align*}
Here, $\sigma_1(\cdot)$ returns the largest singular value of the matrix in its argument.

Note that $L_{\| \cdot \|} [\bf]$ is finite since $\bf$ is Lipschitz continuous.
In this paper, we also assume that $\bf$ is continuously differentiable. In this case, 
$L_{\| \cdot \|} [\bf]= \sup_{\bz} \| D \bf(\bz) \|_{op},$ where $D \bf$ is the Jacobian of $\bf$ and $\|\cdot\|_{op}$ is the operator norm induced by  $\|\cdot\|$. However, for simplicity, we will omit the subscript $_{op}$.

One may use $\tilde{r}_{\| \cdot \|} [\bf] := \log \left( L_{\| \cdot \|} [\bf] \right)$ as an alternative definition of reactivity, which we refer to as ``logarithmic reactivity", 
and obtain the same results that we discuss in what follows. 

Consider a nonlinear and time-varying discrete system
\begin{equation} \label{eq:nonlineardis}
    \bx_{k+1} = \bf_k (\bx_k), \quad k = 0, 1, 2, \hdots, 
\end{equation}
where $\bf_k$s are Lipschitz continuous.
Here, $\bx_k \in \mathbb{R}^m$ is the state of the system at time $k$, and $\bf_k : \mathbb{R}^m \rightarrow \mathbb{R}^m$ is the dynamical function at time $k$.
Let $\bx_k$, $k = 0,  1, 2, \hdots$, denote an orbit of \eqref{eq:nonlineardis} initiated at $\bx_0$.
In what follows, we define the reactivity of \eqref{eq:nonlineardis}. 

\begin{definition} \normalfont{\textbf{(Invariant set)}}
A subset $\mathcal{C} \subseteq \mathbb{R}^m$ is invariant under $\{\bf_k\}$ if $\bx_0\in\mathcal{C}$ implies $\bx_k\in\mathcal{C}$, $\forall k\geq1$.
\end{definition}

\begin{definition} \normalfont{\textbf{($N-$step mean reactivity, sup-mean reactivity, \& lim-mean reactivity of \eqref{eq:nonlineardis}) }}
For a given norm $\|\cdot\|$, we define the $N-$step mean reactivity of  \eqref{eq:nonlineardis} to be the average of the reactivity of $\bf_0, \bf_1,\ldots,\bf_{N-1}$:
   $ \meanR  := 
    \langle r_{\| \cdot \|} [\bf_k]\rangle_{k=0}^{N-1} = \frac{1}{N}\sum_{k=0}^{N-1} r_{\| \cdot \|} [\bf_k],$
the sup-mean reactivity of \eqref{eq:nonlineardis} to be 
    $\barR := \sup_{N} R_{\| \cdot \|}^{(N)},$
and the lim-mean reactivity of \eqref{eq:nonlineardis} to be 
   $ \limR := \lim_{N\to \infty} R_{\| \cdot \|}^{(N)}.$
\end{definition}

Conventionally, a contractive system is thought of as a system with consistent convergence of any of its two trajectories.
However, the reactivity of a system is usually perceived as an indicator of whether the norm of the state vector increases transiently. 
In what follows, we provide a definition of reactive systems that is consistent with current literature on reactivity and contraction theory.

\begin{definition}
 \normalfont{\textbf{(Reactive and Contractive Systems)}}
A system is called sup-mean reactive (resp. sup-mean contractive) on an invariant set $\mathcal{C} \subseteq \mathbb{R}^m$ if there exists a norm $\|\cdot\|$ such that $\barR >0$ (resp. $\barR<0$).
 \end{definition} 

We denote the reactivity of a map $\bf$ 
by the lowercase $r$, the the sup-mean reactivity of 
a discrete dynamical system $\bx_{k+1} = \bf_k (\bx_k)$ by $\barR$, and the lim-mean reactivity of $\bx_{k+1} = \bf_k (\bx_k)$ by $\limR$. Note that, by definition, these notions coincide for time-invariant maps, i.e., for $\bx_{k+1}=\bf(\bx_k)$, $\meanR= \barR =\limR= r[\bf]$. 

For time-invariant maps, our definitions of sup-mean reactivity and sup-mean contractivity are consistent with those for reactivity and contractivity in the literature. However, for time-varying maps, classical reactivity and contractivity are typically defined based on the sign of $R:=\sup_k r_{\|\cdot\|}[\bf_k]$. In this work, we relax this condition: rather than requiring $r_{\|\cdot\|}[\bf_k] >0$  (or $< 0$) for all $k$, we allow the average over $N$ steps to remain positive (or negative). As a result, any reactive or contractive system is also sup-mean reactive or sup-mean contractive, but not necessarily the reverse.
We demonstrate the conservativeness of $r$ over $\barR$ in Sec.~\ref{sec:sync} via an example.

\section{Non-reactivity of the $p$-iteration system}\label{sec:p-th interation}

We begin by stating results on sup-mean contractivity for maps. For time-invariant maps, these results are similar to the contractivity results found in the literature; see \cite{Granas2003,bullo2022contraction}. However, for time-varying maps, these results are novel, as they involve the averages of $r[\bf_k]$ over $N$ steps rather than the point-wise values of $r[\bf_k]$.

\begin{proposition} \label{prop:contractivity}
\normalfont{\textbf{(Sup-Mean Contractive systems)}}
Assume that 
$\mathcal{C}\subset\mathbb{R}^{m}$ is open, convex, and invariant under $\bf_k$.  
Also, assume that there exists a norm $\|\cdot\|$ such that \eqref{eq:nonlineardis} is  sup-mean non-reactive on  $\mathcal{C}$, i.e., $0\leq\lambda:=\barR+1 <1$. 
Then, for any two orbits $\bx=(\bx_0,\bx_1,\hdots)$ and $\by=(\by_0,\by_1,\hdots)$ with $\bx_0\neq\by_0\in\mathcal{C}$, the following inequality holds for any $k\geq1$:
\begin{align}\label{eq:contractivity}
\|\bx_{k} -\by_{k}\| \leq \lambda^k \|\bx_{0} -\by_{0}\|. 
\end{align}
That is, any two orbits converge to each other. 
In addition,  \eqref{eq:nonlineardis} admits at most one fixed point. Moreover, if the map is time independent,  $\bf_k=\bf, \forall k$, and $\mathcal{C}$ is compact, then \eqref{eq:nonlineardis} admits a unique asymptotically stable fixed point in $\mathcal{C}$.
\end{proposition}

\begin{pf}
By Definition~\ref{def:Reactivity1}, we have:
\begin{align*}
\|\bx_{k} -\by_{k}\| 
&= \|\bf_{k-1}(\bx_{k-1}) -\bf_{k-1}(\by_{k-1})\| \\
&\leq (r_{\|\cdot\|} [\bf_{k-1}]+1) \|\bx_{k-1} -\by_{k-1}\|\\
&\qquad\vdots\\
&\leq (r_{\|\cdot\|} [\bf_{k-1}]+1) \cdots(r_{\|\cdot\|} [\bf_0]+1)\|\bx_{0} -\by_{0}\|\\
&\leq \left(\frac{1}{k}\sum_{i=0}^{k-1} r_{\| \cdot \|} [\bf_i] +1\right)^k\;\|\bx_{0} -\by_{0}\|\\
&\leq \lambda^k\|\bx_{0} -\by_{0}\|.
\end{align*}
The second last inequality holds because the arithmetic mean of non-negative real numbers is greater than or equal to their geometric mean.

Next, we show that there exists at most one fixed point. Suppose, for the sake of contradiction, that there are two distinct fixed points, $\bz^* \neq \bar{\bz}$. Then for any $k \geq 1$, $\bf_k(\bar{\bz}) = \bar{\bz}$ and $\bf_k(\bz^*) = \bz^*$. In particular, for $k = 1$,
$
\|\bar{\bz} - \bz^*\| = \|\bf_1(\bar{\bz}) - \bf_1(\bz^*)\| \leq \lambda \|\bar{\bz} - \bz^*\| < \|\bar{\bz} - \bz^*\|.
$
The last inequality follows from $\lambda < 1$. This is a contradiction. Therefore, there cannot be more than one fixed point.

Finally, we assume that $\bf_k = \bf$, and by the Banach Contraction Theorem \cite[Theorem 1.5]{bullo2022contraction}, we show that there exists a unique asymptotically stable fixed point.
The idea is to show that any orbit $\bz_0, \bz_1, \ldots$ is a Cauchy sequence, and therefore converges to a point $\bz^*$ in the compact set $\mathcal{C}$. For any $m > n \geq 1$, we have
\begin{align*}
  \|\bz_{m} - \bz_{n}\| &\leq \|\bz_{m} - \bz_{m-1}\| + \cdots + \|\bz_{n+1} - \bz_{n}\| \\
  &\leq (\lambda^{m-1} + \cdots + \lambda^{n}) \|\bz_{1} - \bz_{0}\| \\
  &\leq \frac{\lambda^{n}}{1 - \lambda} \|\bz_{1} - \bz_{0}\|.
\end{align*}
Since $0\leq\lambda<1$, for any $\epsilon > 0$, there exists $N \geq 1$ such that for all $m > n \geq N$, we have $\|\bz_{m} - \bz_{n}\| < \epsilon$. Thus, the sequence $\{\bz_k\}$ is Cauchy.
\qed 
\end{pf}

Note that for contractive time-dependent maps, the existence of fixed points is not guaranteed. For example, $\bf_k = (-1)^k$ is a contractive map with $\lambda = 0$. However, it does not admit any fixed point; instead, it admits a periodic solution of period 2.

\begin{remark}
\label{remark:stability}
Proposition~\ref{prop:contractivity} addresses contractivity, meaning that over $k$ steps, \eqref{eq:contractivity} is satisfied. However, if the goal is to prove asymptotic stability, i.e., 
$\lim_{k \to \infty} \|\bx_k - \by_k\| = 0, $
the condition in Proposition~\ref{prop:contractivity} can be relaxed to assume that $0\leq\limR+1 <1$.

The reason is that by \eqref{eq:contractivity}, for $k\geq1$ we have:
$
\|\bx_{k} -\by_{k}\| 
 \leq \left(\frac{1}{k}\sum_{i=0}^{k-1} r_{\| \cdot \|} [\bf_i] +1\right)^k\;\|\bx_{0} -\by_{0}\|. 
$
Since $\lim_{k\to\infty}\frac{1}{k}\sum_{i=0}^{k-1} r_{\| \cdot \|} [\bf_i] +1 = \limR+1<1$, for any $0<\epsilon<-\limR$, there exists $N\geq1$ such that for $k\geq N$, 
$\frac{1}{k}\sum_{i=0}^{k-1} r_{\| \cdot \|} [\bf_i] +1 < \limR+1 +\epsilon<1$, and hence,
$\left(\frac{1}{k}\sum_{i=0}^{k-1} r_{\| \cdot \|} [\bf_i] +1\right)^k < (\limR+1 +\epsilon)^k$. 
Therefore, 

\begin{align*}
\begin{split}
 \lim_{k\to\infty}\|\bx_{k} -\by_{k}\|
 &\leq \lim_{k\to\infty}\left(\frac{1}{k}\sum_{i=0}^{k-1} r_{\| \cdot \|} [\bf_i] +1\right)^k\;\|\bx_{0} -\by_{0}\|\\
 &\leq \lim_{k\to\infty}(\limR+1 +\epsilon)^k
    = 0.
\end{split}
\end{align*}
The last inequality holds since $0\leq\limR+1 +\epsilon<1$.
\end{remark}

In what follows, we prove our main result, that is, if the first iteration system is not sup-mean contractive but there exists $p>1$ such that the $p$-th iteration of that system is sup-mean contractive, then the system admits an asymptotically stable attractor.
First, we define the $p$-iteration system for an integer $p \geq 1$:
\begin{definition}\label{def:p-iteration}
    The $p$-iteration system corresponding to the dynamical system \eqref{eq:nonlineardis} is defined as
    \begin{align*}
        \bx_{x + p} = \bf_{k}^{(p)}(\bx_k) &:= 
         \bf_{k+p-1} \circ \bf_{k+p-2}\circ \cdots\circ \bf_{k} (\bx_k)\\
        & =\bf_{k+p-1} (\bf_{k+p-2}( \cdots ( \bf_{k} (\bx_k)) ) ).
    \end{align*}
    Note that the $1$-iteration system $\bf_{k}^{(1)}(\bx_k) = \bf_k(\bx_k)$. For $2$-iteration system $\bf_{k}^{(2)}(\bx_k) = \bf_{k+1}(\bf_k(\bx_k))$.
    
\end{definition}
It is known \cite{nazerian2023reactability} that if the reactivity of a linear time-invariant discrete dynamical system is negative, i.e., if the system is \textit{sup-mean non-reactive} (or \textit{sup-mean contractive}) then the system is asymptotically stable. Our main goal is to prove that if a system is not sup-mean contractive but its $p$-th iteration is sup-mean contractive, it still admits an asymptotically stable attractor (e.g., fixed points, periodic solutions, etc.).

\subsection{Time-varying linear maps}
\label{sec:time-variant-linear-maps}
 Consider a linear time-varying discrete-time system (which is for example the case of a linearized discrete-time system about a given trajectory)
\begin{equation} \label{eq:discrete}
    \bx_{k+1} = A_k \bx_k,
\end{equation}
for which, by Definition \ref{def:p-iteration}, the $p$-iteration is
\begin{equation} \label{eq:discretep}
    \bx_{k+p} = \mathcal{A}_k^{(p)} \bx_k
\end{equation}
with
\begin{equation} \label{eq:piter}
    \mathcal{A}_k^{(p)} = A_{k+p-1} A_{k+p-2} \cdots A_{k+1} A_k.
\end{equation}

\begin{remark}
    It is known that 
    \begin{align*}
        \| \mathcal{A}_k^{(p)} \| = \|  A_{k+p-1} \cdots & A_{k+1} A_k \| 
        \\ & \leq \|  A_{k+p-1} \| \cdots \| A_{k+1} \| \| A_k \|.
    \end{align*}
    If the terms on the right-hand side are all less than one, it follows that the $\| \mathcal{A}_k^{(p)} \| < 1$. However, the converse is not true.  Therefore, it is possible that $\| A_{i} \| > 1$ for some $k \leq i \leq k+p-1$, but $\| \mathcal{A}_k^{(p)} \| < 1$.
    This property provides the motivation to study $\| \mathcal{A}_k^{(p)} \|$ rather than $\| A_{k} \|$.
\end{remark}

\begin{proposition} \label{prop:time_variant_linear}
    Consider \eqref{eq:discrete} and \eqref{eq:discretep} and assume that $\exists p \geq 1$ such that \eqref{eq:discretep} is sup-mean contractive, i.e., $\barR = \sup_N \frac{1}{N}\sum_{k=0}^{N-1} r_{\| \cdot \|} [\mathcal{A}_k^{(p)}] < 0$ for some norm $\|\cdot\|$. Then \eqref{eq:discrete} is asymptotically stable.
\end{proposition}
\begin{pf}
    Let $\{ \bx_l \}_{l=0}^\infty$ be the orbit of \eqref{eq:discrete} with initial condition $\bx_0$. 
    We want to show that $\{ \bx_l \}_{l=0}^\infty$ converges to zero, i.e., $\forall \epsilon > 0, \ \exists L > 0$ such that for $l\geq L$, $\| \bx_l \| < \epsilon$.
Consider \eqref{eq:discretep} with initial conditions $\{\bx_{0}, \ldots,\bx_{p-1}\}$. Then 
\begin{align} \label{eq:proof}
        \{ \bx_l \}_{l=0}^\infty = \cup_{k = 0}^{p-1} \{ \bx_{k+jp} \}_{j = 0}^\infty,
    \end{align}
    \noindent where $\{ \bx_{k+jp} \}_{j = 0}^\infty$ is the trajectory of \eqref{eq:discretep} with initial condition $\bx_{k}$.
    Since \eqref{eq:discretep} is contractive and admits only one fixed point at the origin, based on Proposition \ref{prop:contractivity}, the fixed point is asymptomatically stable. That is, $\forall \epsilon > 0, \ \exists M_k > 0$, such that for $j \geq M_k$, $\| \bx_{k+jp}\| < \epsilon$.
    Let $M = \max \{ M_0, M_2, \hdots, M_{p-1} \}$. 
    Then $\forall j \geq M$ and $k = 0, \hdots, p-1$, $\| \bx_{k+jp} \| < \epsilon$.
    Using \eqref{eq:proof}, we conclude that $\| \bx_l \| < \epsilon$ for $l \geq L:=pM$.
\qed 
\end{pf}

\textbf{Example 1.} 
Consider the system $\bx_{k+1} = A_k \bx_k$ with
$
    A_k = \begin{bmatrix}
        0.5 & \lambda^k \\
        0 & 0.5
    \end{bmatrix}.
$
It is known \cite[Example 4]{Zhou2017Asymptotic} that this system is asymptotically stable if and only if $| \lambda | < 2$.
For $|\lambda|\leq1$ and $k\geq0$, 
$r_{1,Q}[A_k] = 0.5 + |\lambda|^k/3 -1<0$, where $r_{1,Q}$ is the reactivity induced by $Q-$weighted $L_1$ norm with $Q=\mbox{diag}\{1,3\}$. This implies asymptotic stability.

\begin{figure}
    \centering
    \includegraphics[width=0.9\linewidth]{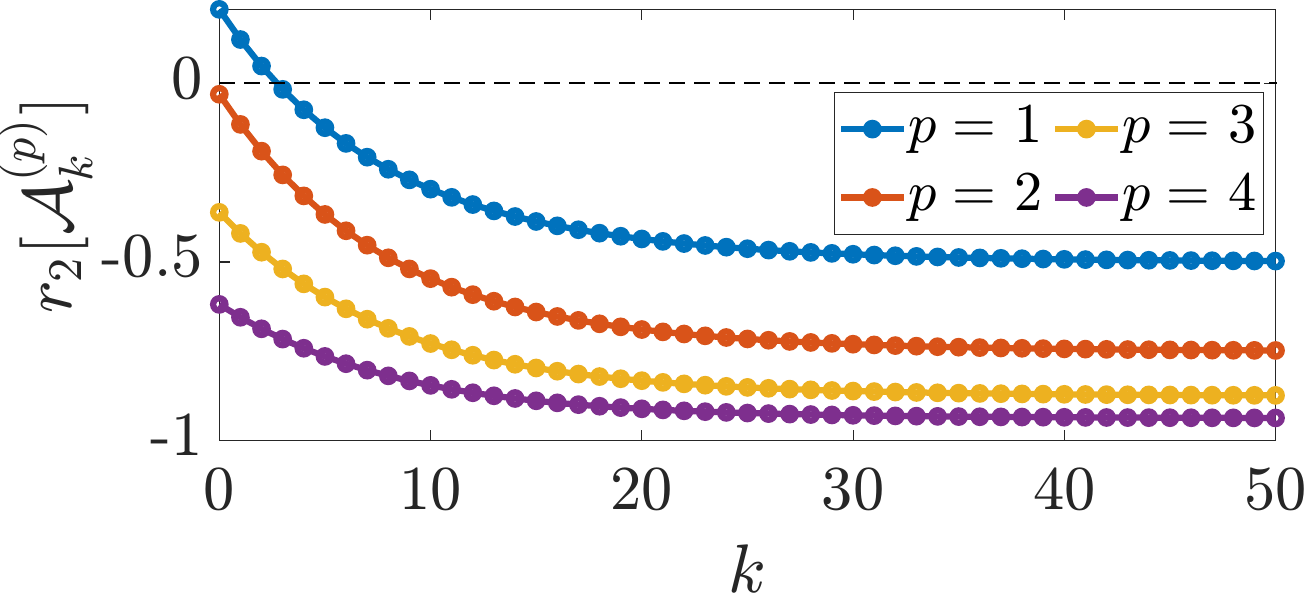}
    \caption{The reactivity of the $p$-iteration system. For $p= 1$, $r_2[\mathcal{A}_k^{(p)}] > 0$ for $k = 0,1, 2$. Here, $\lambda = 0.9$ in Example 1. 
    }
\label{fig:example2}
\end{figure}

The $p$-iteration system is characterized by the matrix
\begin{align} \label{eq:Ap}
\begin{split}
    \mathcal{A}_k^{(p)} & = \begin{bmatrix}
        (0.5)^p & (0.5)^{p-1}\lambda^k \left( 1 + \lambda  + \cdots \lambda^{p-1} \right) \\
        0 & (0.5)^p
    \end{bmatrix} 
    \end{split}
    \end{align}
Here, we set $\lambda = 0.9$.
The reactivity of $\bx_{k+1} = A_k \bx_k$ is $r_2[A_k] = \sigma_1(A_k) - 1 > 0, k = 0, 1, 2$. 
However, for $p = 2$, the reactivity $r_2 (\mathcal{A}_k^{(p)}) < 0, \forall k \geq 0$.
Figure\,\ref{fig:example2} shows the values of $r_2[\mathcal{A}_k^{(p)}]$ for $p = 1, 2, 3, 4$ and $0\leq k \leq 50$.
As can be seen in this example, the study of the $p=2$ iteration system is sufficient to assess the stability of the system since if $r_2[\mathcal{A}_k^{(2)}] < 0$, $\forall k$, then $\bar{R}_2 = \sup_N \frac{1}{N}\sum_{k=0}^{N-1} r_{2} [\mathcal{A}_k^{(2)}] < 0$.

Here, we calculate $r_{\infty}[\mathcal{A}_k^{(p)}]$ for $p \rightarrow \infty$ and $k \rightarrow \infty$:
From \eqref{eq:Ap}, we know
$
    r_{\infty}[\mathcal{A}_k^{(p)}] = (0.5)^{p} + (0.5)^{p-1} \left| \lambda^k \frac{1-\lambda^p}{1-\lambda} \right|.
$
If $p \rightarrow \infty$ and $k \rightarrow \infty$, then
$
    r_{\infty}[\mathcal{A}_\infty^{(\infty)}] = (0.5)^{\infty} \left|  \frac{\lambda^\infty}{1-\lambda} \right|.
$
Thus, $r_{\infty}[\mathcal{A}_\infty^{(\infty)}] < 1$ if $0.5 | \lambda | < 1$ which means $| \lambda | < 2$ as confirmed in \cite[Example 4]{Zhou2017Asymptotic}.

\subsection{Time-invariant nonlinear maps}
\label{sec:nonlinear-maps}

Consider a nonlinear time-invariant map $\bf$. In the following proposition, we show that if the $p-$iteration of $\bf$ is contractive on an invariant subset $\mathcal{C}$, then $\bf$ admits a unique asymptotically stable fixed point on $\mathcal{C}$.
References \cite[Theorem 1]{bessaga1959converse} and \cite[Section 1.6 (A.1)]{Granas2003} also exploited the notion of $p$-iteration system to state contraction conditions closely related to Proposition~\ref{prop:nonlinear}. 
However, those references do not provide proof.

\begin{proposition} \label{prop:nonlinear}
Assume that $\mathcal{C}$ is a convex and compact subset of $\mathbb{R}^m$ and it is invariant under $\bf$. 
Assume that there exist $p\in\mathbb{N}$ and a norm $\|\cdot\|$ such that $\bx_{k+p}= \bf^{(p)}(\bx_k)$ is sup-mean contractive in $\mathcal{C}$, i.e., $\bar R_{\|\cdot\|} = r[\bf^{(p)}]<0$. 
Then, $\bx_{k+1}= \bf(\bx_k)$ has a unique globally asymptotically stable fixed point in $\mathcal{C}$.
\end{proposition}

\begin{pf}
Existence of a unique fixed point: Since $\mathcal{C}$ is invariant for $\bf$, we have
$\bx\in\mathcal{C} \Rightarrow
\bf(\bx)\in\mathcal{C}
\Rightarrow
\cdots
\Rightarrow
\bf^{(p)}(\bx)\in\mathcal{C}.
$
Therefore, $\mathcal{C}$ is invariant for $\bf^{(p)}$.
Since $\bf^{(p)}$ is sup-mean contractive in compact set $\mathcal{C}$, by Proposition~\ref{prop:contractivity}, it admits a unique asymptotically stable fixed point $\bx^*\in\mathcal{C}$, i.e., $\bf^{(p)} (\bx^*) =\bx^*$. We prove that $\bx^*$ is the unique fixed point of $\bf$. Let $\bf(\bx^*)=\by^*$. Since $\bx^*\in\mathcal{C}$,  $\by^*\in\mathcal{C}$, 
$
    \bf^{(p)} (\by^*)=\bf^{(p)} (\bf(\bx^*))
                     = \bf(\bf^{(p)}(\bx^*))
                     =\bf(\bx^*)
                     =\by^*. 
$
Since $\by^*\in\mathcal{C}$ and $\bx^*$ is the unique fixed point of $\bf^{(p)}$ in $\mathcal{C}$, we conclude that $\by^*=\bx^*$, and hence, $\bx^*$ is a fixed point of $\bf$. Note that this is a unique fixed point because if $\bf$ has another fixed point $\bz^*$, then $\bz^*$ is also a fixed point of $\bf^{(p)}$ which contradicts the uniqueness of $\bx^*$.   

The asymptotic stability of the unique fixed point is similar to that established in Proposition~\ref{prop:time_variant_linear}, so we omit the details.
Let $\{\bx_k\}_{k=0}^\infty$ denote the orbit of the map $\bx_{k+1} = \bf(\bx_k)$ with initial condition $\bx_0 \in \mathcal{C}$. We will show that $\bx_k \to \bx^*$ as $k \to \infty$.
Consider $p$ orbits of the $p$-iteration map $\bx_{l+p} = \bf^{(p)}(\bx_l)$, each with initial conditions $\bx_0, \ldots, \bx_{p-1}$. Then, the orbit $\{\bx_k\}_{k=0}^\infty$ can be viewed as the union of these $p$ orbits. Since each orbit converges to $\bx^*$ asymptotically, their union also converges to $\bx^*$ asymptotically.\qed
\end{pf}

Proposition~\ref{prop:nonlinear} states that if, for some $p>1$, the $p$-iteration of $\bf$ is contractive on an invariant set $\mathcal{C}$, then $\bf$ is convergent in $\mathcal{C}$.
 In Proposition~\ref{prop:limit_cycles_time_independent}, we show that $\bf^{(p)}$ is not necessary to be contractive on the whole set $\mathcal{C}$. In fact, if it is contractive in a subset of $\mathcal{C}$ that is invariant under the dynamics of $\bf^{(p)}$, then we can still derive convergence results for $\bf$.

\begin{proposition}\label{prop:limit_cycles_time_independent}
Assume that $\mathcal{C}_p\subset\mathcal{C}_1$ are compact subsets of $\mathbb{R}^m$ and invariant under the orbits of $\bf^{(p)}$ and $\bf$, respectively. In addition, $\bf^{(p)}$ is sup-mean contractive on $\mathcal{C}_p$. That is, $\bx\in\mathcal{C}_p\Rightarrow \bf^{(p)}(\bx)\in\mathcal{C}_p$ and $\bar R_{\|\cdot\|} = r[\bf^{(p)}]<0$, for some norm $\|\cdot\|$. Then $\bf$ admits either a unique stable fixed point in $\mathcal{C}_p$ or a stable periodic solution in $\mathcal{C}_1$ of period $d$ where $d$ is a divisor of $p$. 
\end{proposition}

\begin{pf} 
Since $\bf^{(p)}$ is sup-mean contractive on $\mathcal{C}_p$, by Proposition~\ref{prop:contractivity}, it admits a unique fixed point $\bx^*\in\mathcal{C}_p$, i.e., $\bf^{(p)}(\bx^*) = \bx^*$. Now, let 
$\bz_1^* := \bf(\bx^*), 
\bz_2^*: =  \bf(\bz_1^*) =\bf^{(2)}(\bx^*), \ldots,
\bz_p^*: =  \bf(\bz_{p-1}^*) =\bf^{(p)}(\bx^*) = \bx^*$. 
Since $\bx^* = \bz_p^*$, we can conclude that 
$\bx^* = \bz_d^*$, for some divisor $d$ of $p$. In this case, if $d\neq1$, then $\bf$ admits a periodic solution of period $d$ and 
if $d=1$, then $\bx^*$ is the unique fixed point for $\bf$. 
Note that if $\bf(\by^*) = \by^*$ for some $\by^* \neq \bx^*$ in $\mathcal{C}_p$, then $\bf^{(p)}(\by^*) = \by^*$, which contradicts the uniqueness of the fixed point of $\bf^{(p)}$ in $\mathcal{C}_p$. 
\qed
\end{pf}

\textbf{Example 2.} (Time-invariant Logistic Map)
The following example demonstrates the results provided in Proposition \ref{prop:limit_cycles_time_independent}.
The Logistic map $f:[0,1]\to[0,1]$ is governed by 
$
    x_{k+1} = f(x_k) = \alpha x_k (1 - x_k),
$
where $\alpha$ is a parameter. 
\ZA{For $\alpha$ in the interval $[1, 4]$, $x_n$ remains in $[0, 1]$.}
The map admits two fixed points: 
$\bx^*_1=0$, which is always unstable, and $\bx^*_2 = (\alpha - 1) / \alpha$, which is stable when $\alpha\in [1,3]$ and loses stability, giving rise to a periodic solution, when $\alpha\in (3,4]$. 

\begin{figure*}
    \centering
    \text{(a): interval width} \hspace{0.18\linewidth} \text{(b): $p=2$} \hspace{0.2\linewidth} \text{(c): $p=4$} \\
    \includegraphics[width=0.29\linewidth]{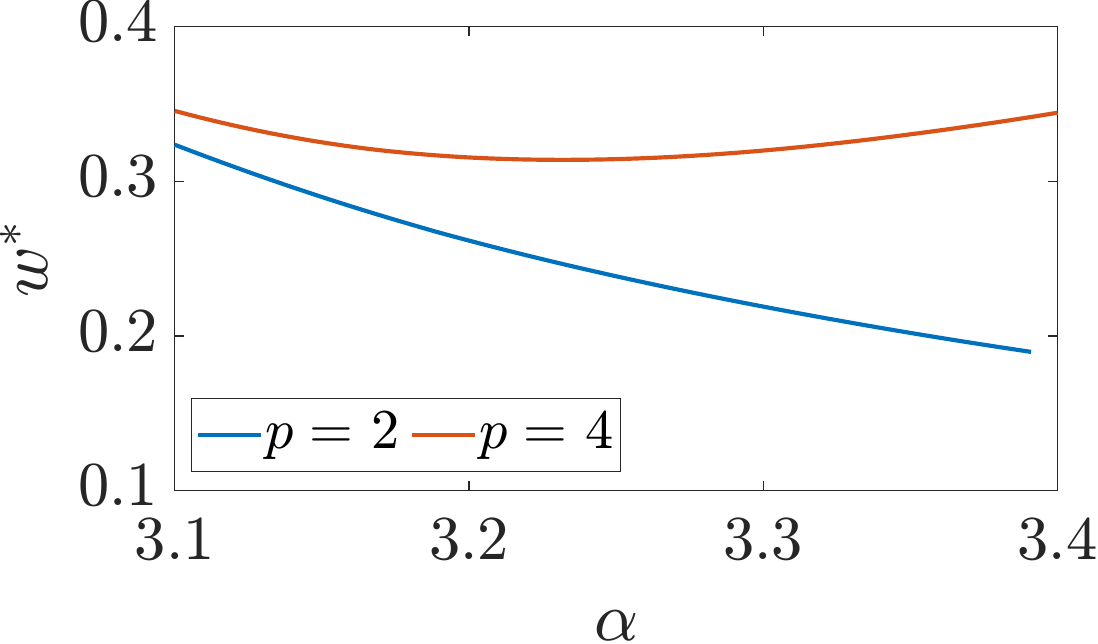}
    \hspace{0.04\linewidth}
    \includegraphics[width=0.29\linewidth]{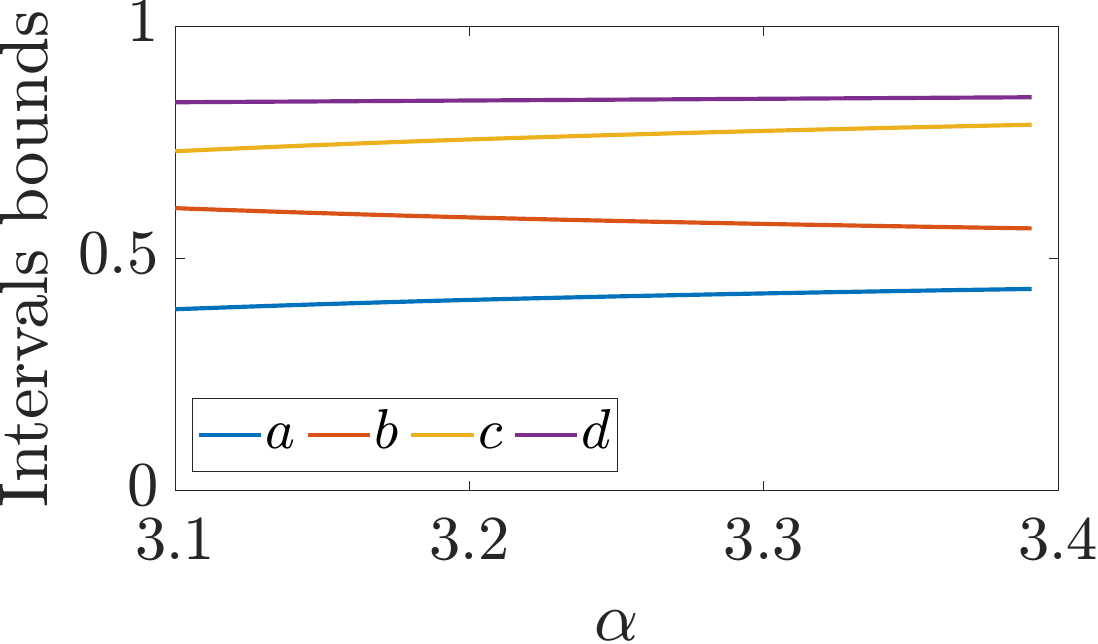} 
    \hspace{0.04\linewidth}
    \includegraphics[width=0.29\linewidth]{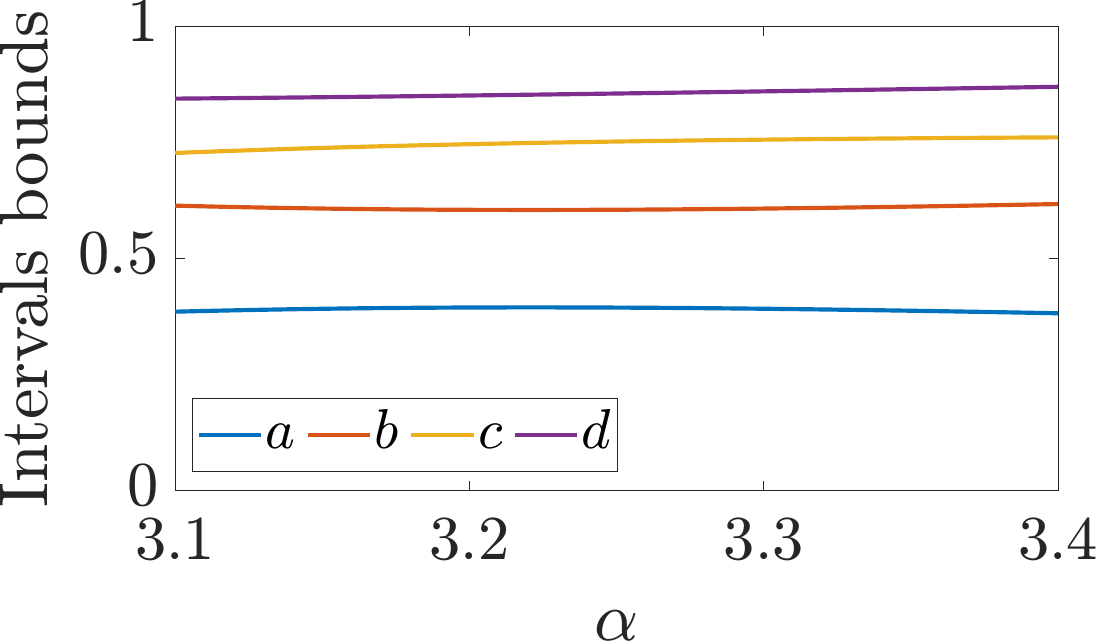}
    
    \caption{(a) The maximum value $w^*$ from \eqref{eq:optperiodic}; (b,c) the optimal values of $a, b, c, d$ from \eqref{eq:optperiodic} as functions of the parameter $\alpha$ of the time-invariant Logistic map in Example 2, for iteration numbers $p = 2$ and $p = 4$, respectively.
    The gaps within the curves correspond to the infeasibility of the optimization for $\alpha$ close to $3.4$.}
    \label{fig:periodic}
\end{figure*}

We are interested in analyzing the stability of the periodic Logistic map.
We seek the widest intervals $\mathcal{C}_1 = [a, \ b] \cup [c, \ d]$ and $\mathcal{C}_p = [a, \ b]$ such that $c > b$ and $f$ satisfies $f([a, \ b]) \subseteq [c, \ d]$ and $f([c, \ d]) \subseteq [a, \ b]$. \ZA{This implies that $f$ is invariant in $\mathcal{C}_1$ and $f^{(p)}$ is invariant in $\mathcal{C}_p$ when $p$ is even. }
We also require $\| Df^{(p)}\| \leq 1,$ on $\mathcal{C}_p$.
To find $\mathcal{C}_1$ and $\mathcal{C}_p$ for a given parameter $\alpha$ and iteration $p$, we solve: 
\begin{align} \label{eq:optperiodic}
\begin{split}
    w^* = \max_{a,b,c,d} \quad & b -a + d - c \\
    \text{s.\,t.} \quad & 0 \leq a \leq b \leq (\alpha - 1) / \alpha - \epsilon,  \\
    & (\alpha - 1) / \alpha + \epsilon \leq c \leq d \leq 1, \\
    & f: [a,\, b] \to [c,\, d], \\
    & f : [c,\, d] \to [a,\, b], \\
    & \|Df^{(p)}\| \leq 1, \quad\text{on}\; \mathcal{C}_p,
\end{split}
\end{align}
where $\epsilon$ is a small positive number to satisfy the strict inequality condition $c > b$. Here we set $\epsilon = 10^{-6}$.
We vary $3.1 \leq \alpha \leq 3.4$ and set $p = 2, 4$. 
\color{black}
Figure \ref{fig:periodic} (left panel) shows $w^*$ after solving \eqref{eq:optperiodic} as the parameter $\alpha$ and the iteration $p$ are varied. 
It is clear that by the use of higher $p$, a wider interval of stability based on Proposition \ref{prop:limit_cycles_time_independent} for the periodic orbit of the Logistic map is predicted. 
Figure \ref{fig:periodic} (center and right panels) shows the optimal bounds $a,b,c$, and $d$ of the intervals $\mathcal{C}_1$ and $\mathcal{C}_p$ after solving \eqref{eq:optperiodic} for $p=2$ and $p=4$.

\subsection{Time-Varying Nonlinear Maps}

In what follows, we generalize the result of Proposition~\ref{prop:limit_cycles_time_independent} to a special case of time-varying discrete systems. 

\begin{proposition}\label{prop:limit_cycles_time_dependent}
Let $\mathcal{C}_1$ and $\mathcal{C}_2$ be two convex subsets of $\mathbb{R}^m$ and
$\bg:\mathcal{C}_1\to\mathcal{C}_2$ and  $\bh:\mathcal{C}_2\to\mathcal{C}_1$ be Lipschitz continuous functions. 
Assume that 
$\bh\circ\bg:\mathcal{C}_1\to\mathcal{C}_1$ and $\bg\circ\bh:\mathcal{C}_2\to\mathcal{C}_2$ are sup-mean contractive.   
Now, for $k\geq0$ define the time-varying map $\bf_k$ as follows: $\bf_{2k}=\bg$ and $\bf_{2k+1}=\bh$. Then, the orbit of $\bz_{k+1} = \bf_k(\bz_k)$ with initial value $\bz_0\in\mathcal{C}_1$ converges either to a limit cycle of period 2 or a fixed point. 
\end{proposition}
\begin{pf}
By definition 
$\bh\circ\bg$ and $\bg\circ\bh$ are invariant in $\mathcal{C}_1$ and $\mathcal{C}_2$, respectively. Since $\bh\circ\bg$ and $\bg\circ\bh$ are also contractive in $\mathcal{C}_1\cup\mathcal{C}_2$, by Proposition~\ref{prop:contractivity}, they admit unique asymptotically stable fixed points in the corresponding sets. That is,  $\exists\;\bx^*\in \mathcal{C}_1$ such that $\bh(\bg(\bx^*)) = \bx^*$ and $\exists\; \by^*\in \mathcal{C}_2$ such that $\bg(\bh(\by^*)) = \by^*$.

 Let $\{\bz_k\}_{k=0}^{\infty}$ be an orbit of $\{\bf_k\}_{k=0}^{\infty}$ where $\bz_0\in\mathcal{C}_1$. First, we show that $\{\bz_k\}$ is convergent.
By definition, $\by_k:=\bz_{2k}$ is an orbit of $\bh\circ\bg$ which converges to $\by^*$ and 
$\bx_k:=\bz_{2k+1}$ is an orbit of $\bg\circ\bh$ which converges to $\bx^*$.
 Therefore, $\{\bz_k\}$ converges to $\{\bx^*,\by^*\}$. 
 \; Next, we show that $\bg(\bx^*)=\by^*$ and 
$\bh(\by^*)=\bx^*$. This proves the existence of the limit cycle. 
Assume that $\bg(\bx^*)=\bz^*$ for some $\bz^*\in\mathcal{C}_2$. Then (by applying $\bh$ to both sides) $\bx^*=\bh(\bg(\bx^*))=\bh(\bz^*)$, and hence (by applying $\bg$ to both sides), 
$\bz^*=\bg(\bx^*)=\bg(\bh(\bz^*))$. Since $\by^*$ is the unique fixed point of $\bg\circ\bh$, we conclude that $\bz^*=\by^*$. Similarly, we can show that $\bh (\by^*) = \bx^*$. 

Then, $\{\bz_k\}$ converges to a fixed point if $\bx^*=\by^*$ and a limit cycle of period 2 if $\bx^*\neq\by^*$.
\qed
\end{pf}

\begin{remark}
Proposition~\ref{prop:limit_cycles_time_dependent} can be extended to the case where the time-varying map is composed of more than two functions, that is, $\bf_{mk+j} = \prescript{}{j}\bh$ for some functions $\prescript{}{1}\bh, \ldots, \prescript{}{m}\bh$ with $m > 2$. Here, the index $j$ selects from $m>2$ time-independent functions. In this setting, if each composition $\prescript{}{j+1}\bh \circ\prescript{}{j}\bh$ is contractive, then  $\bf_k$ admits a limit cycle of period $d$, where $d\geq1$ is a divisor of $m$. 

Proposition~\ref{prop:limit_cycles_time_dependent} can also be extended to the case of $p$-iteration, where $p$ is even.
 
\end{remark}

\begin{figure*}
    \centering
    \text{(a): interval width} \hspace{0.18\linewidth} \text{(b): $p=2$} \hspace{0.2\linewidth} \text{(c): $p=4$} \\
    \includegraphics[width=0.29\linewidth]{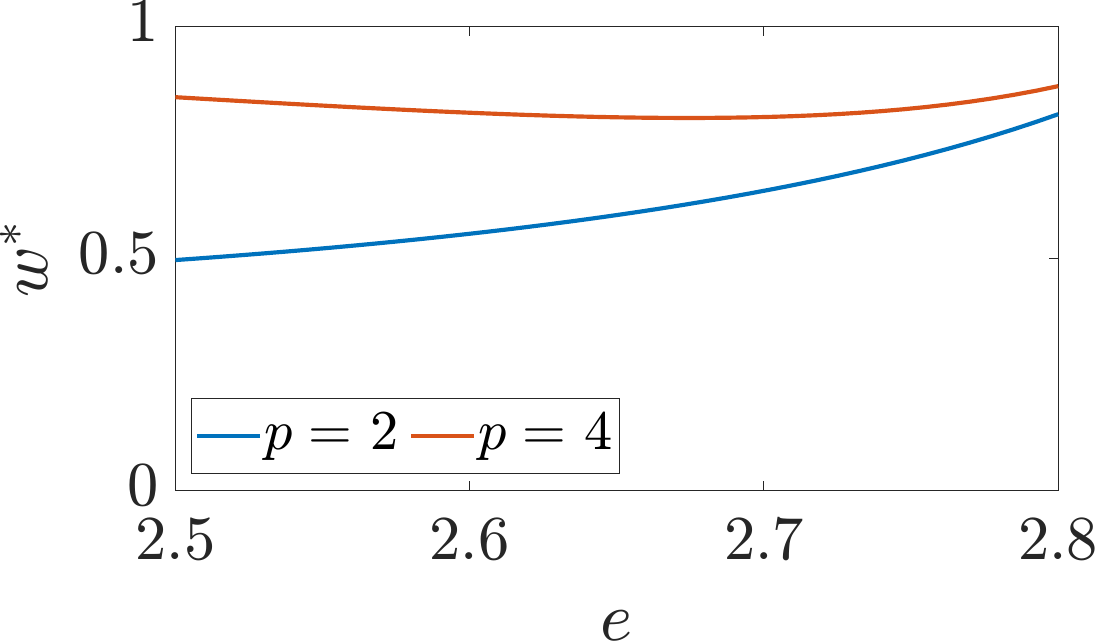}
    \hspace{0.04\linewidth}
    \includegraphics[width=0.29\linewidth]{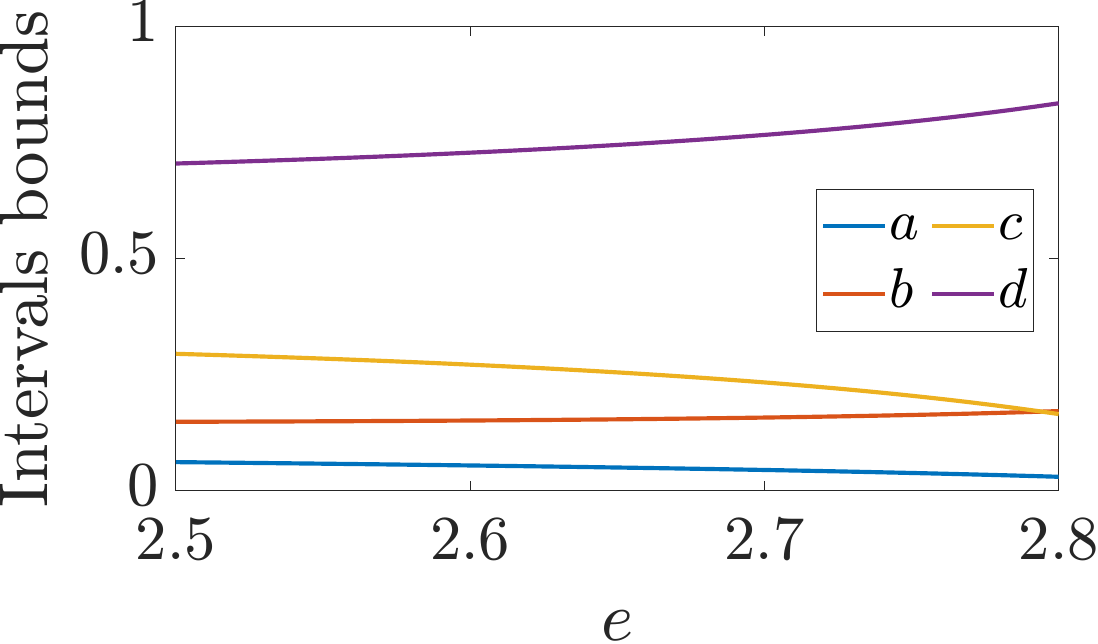} 
    \hspace{0.04\linewidth}
    \includegraphics[width=0.29\linewidth]{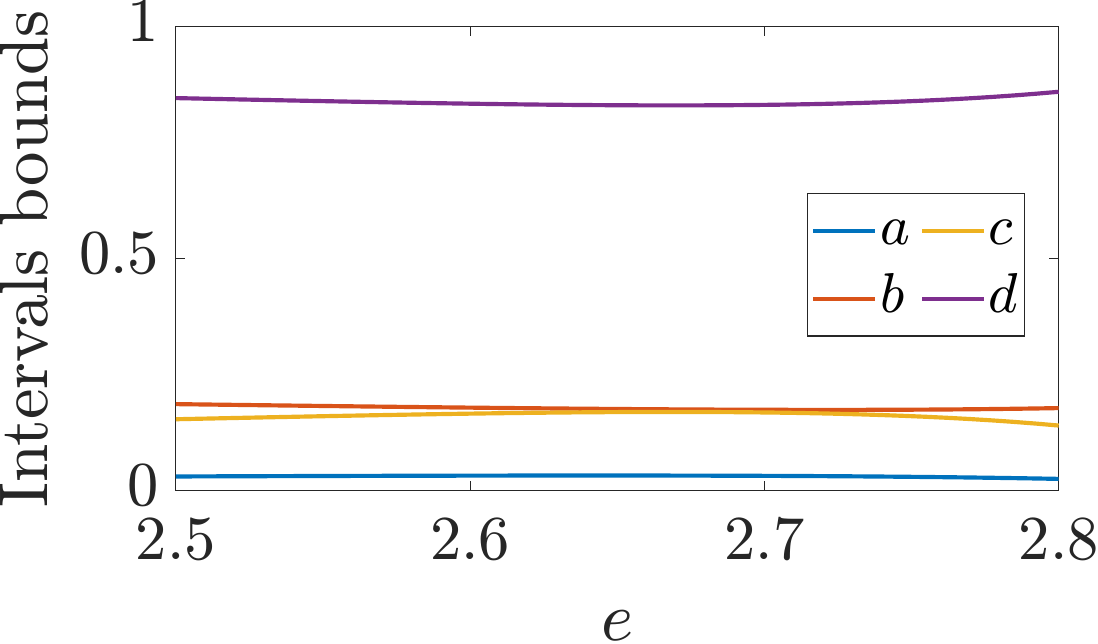}
    
   \caption{(a) The maximum value $w^*$ from \eqref{eq:optperiodic2}; (b,c) the optimal values of $a, b, c, d$ from \eqref{eq:optperiodic2} as functions of the parameter $e$ of the time-varying Logistic map \eqref{eq:logistictv}, for iteration numbers $p = 2$ and $p = 4$, respectively. }
    \label{fig:periodic2}
\end{figure*}

\textbf{Example 3.} (Time-varying Logistic Map) Here, we present an example involving time-varying nonlinear maps. 
Following the notation in Proposition~\ref{prop:limit_cycles_time_dependent}, we set $g(x) = \left( 3.075 + e \right) x (1-x)$, and $h(x) = \left( 3.075 - e \right) x (1-x)$,
and define the time-varying Logistic map as
\begin{equation} \label{eq:logistictv}
    x_{k+1} = f_k(x_k) = \begin{cases}
        g(x_k), \quad & k=0,2,4,\hdots, \\
        h(x_k), \quad & k=1,3,5,\hdots.
    \end{cases}
\end{equation}
The system switches between two values of the time-varying parameter $\alpha_k$ of the Logistic map, i.e.,  $\alpha_k= 3.075 + (-1)^k e$.
For $2.5 \leq e \leq 2.9$, the system has a stable period-2. 
We aim to find the widest pair of intervals $\mathcal{C}_1 = [a, b]$ and $\mathcal{C}_2 = [c, d]$ such that $g: \mathcal{C}_1 \to \mathcal{C}_2$ and $h: \mathcal{C}_2 \to \mathcal{C}_1$.
For $l\in \mathbb{Z}_{\geq0}$ and an even number $p$, we define
\[
    f^{(p)}_{2l} := \underbrace{h\circ g\circ\cdots\circ h\circ g}_{p \text{ factors}}\;:\; \mathcal{C}_1\to \mathcal{C}_1,
    \]
    and
\[
    f^{(p)}_{2l+1} := \underbrace{g\circ h\circ\cdots\circ g\circ h}_{p \text{ factors}}\;:\; \mathcal{C}_2\to \mathcal{C}_2,
\]
and require $\|Df_{2l}^{(p)}\|\leq1$ and $\|Df_{2l+1}^{(p)}\|\leq1$, 
for $p = 2$ and $p=4$.
To find $\mathcal{C}_1$ and $\mathcal{C}_2$ for a given parameter $e$ and the iteration $p$, we solve:
\begin{align} \label{eq:optperiodic2}
\begin{split}
    w^* = \max_{a,b,c,d} \quad & b -a + d - c \\
    \text{s.\,t.} \quad & 0 \leq a \leq b,  \\
    & c \leq d \leq 1, \\
    & g: [a,\, b] \to [c,\, d],\\
    & h: [c,\, d]\to [a,\, b], \\
    & \left\| Df^{(p)}_{2l} \right\| \leq 1, \;\text{ on } [a, \, b], \\
    & \left\| Df^{(p)}_{2l+1}\right\| \leq 1, \;\text{ on } [c,\, d].
\end{split}
\end{align}
We vary $2.5 \leq e \leq 2.8$ and set $p = 2, 4$.
Figure \ref{fig:periodic2} (a) shows $w^*$ after solving \eqref{eq:optperiodic2} as the parameter $e$ and iteration $p$ are varied. 
It is clear that by using a higher $p$, a wider interval of stability is predicted based on Proposition~\ref{prop:limit_cycles_time_dependent} for the periodic orbit of the logistic map. 
Figure~\ref{fig:periodic2} (b-c) shows the optimal bounds $a,b,c$, and $d$ of the intervals $\mathcal{C}_1$ and $\mathcal{C}_2$ after solving \eqref{eq:optperiodic2}.

\section{Connection with Lyapunov Exponent} \label{sec:cle}
Consider a discrete-time map $\bx_{k +1 } = \bf(\bx_k)$ 
where the state vector is $\bx_k \in \mathbb{R}^m$, and $\bf : \mathbb{R}^m \to \mathbb{R}^m$ is a nonlinear differentiable map.
Assume the system is initialized at $\bx_0$.
The time evolution of the infinitesimal perturbations $\by_k$ 
along the orbit $\bx_k$ is
$
    \by_{k+1} = \bD \bf (\bx_k) \by_k.
$
The $p$ step perturbations are $\by_p = \bD \bf^{(p)} (\bx_0) \by_0$ where
$
    \bD \bf^{(p)} (\bx_0) = \bD \bf (\bx_{p-1}) \bD \bf (\bx_{p-2}) \cdots  \bD \bf (\bx_{0}).
$
The notation $\bD \bf^{(p)}$ coincides with the general notation $\mathcal{A}_k^{(p)}$ introduced in \eqref{eq:piter}.
The finite-time Lyapunov exponent \ZA{(FTLE)} is defined for the initial condition $\bx_0$ and the initial orientation of the infinitesimal perturbations $\bu_0 = \by_0 / \| \by_0 \|_2$ as (see \cite[Chapter 4.4]{ott2002chaos})
\begin{align*}
        h(\bx_0, \bu_0) = & \frac{1}{p} \log \left(\frac{\| \by_p \|_2}{\| \by_0 \|_2} \right)  =  \frac{1}{p} \log \| \bD \bf^{(p)} (\bx_0) \bu_0 \|_2 \\
    = & \frac{1}{2p} \log \sqrt{ \bu_0^\top \bD \bf^{(p)} (\bx_0)^\top \bD \bf^{(p)} (\bx_0) \bu_0 }.
\end{align*}
Since the length of the vector $\bx$ is $m$, there are $m$ or less distinct FTLEs. 
The maximum FTLE (MFTLE) corresponds to the largest eigenvalue of the matrix $\bD \bf^{(p)} (\bx_0)^\top \bD \bf^{(p)} (\bx_0)$, or equivalently, the largest singular value of the matrix $\bD \bf^{(p)} (\bx_0)$. 
Hence, 
\begin{equation} \label{eq:mftle}
    MFTLE(\bx_0) =  \frac{1}{2p} \log \sigma_1 (\bD \bf^{(p)} (\bx_0)).
\end{equation}
The reactivity associated with the $p$-iteration induced by the $L_2$ norm is
$
    r_2 \left[\bD \bf^{(p)} (\bx_0) \right] = 
    \sigma_1 \left(\bD \bf^{(p)} (\bx_0) \right) -1.
$
Hence, 
   $ MFTLE(\bx_0)  =  \frac{1}{2p} \log (r_{2}[\bD \bf^{(p)} (\bx_0)] +1),$
which explicitly connects the maximum Lyapunov exponent to the $p$-iteration system.

\begin{remark}
    In the limit of $p \rightarrow \infty$, the maximum FTLE from \eqref{eq:mftle} becomes independent of the initial condition $\bx_0$ and is then referred to as the Maximum Lyapunov Exponent \cite{ott2002chaos}. 
    Similarly, we conclude that in \eqref{eq:piter} as $p \to \infty$, 
$r_2[\mathcal{A}^{(p)}_k] \to \hat{r},$  for some $\hat{r} > 0$ which is independent of $k$, i.e., as long as $p$ is sufficiently large, any choice for $k$ in \eqref{eq:piter} results in approximately the same reactivity of $\mathcal{A}^{(p)}_k$ with $L_2$ norm.
\end{remark}

Given that the $p$-iteration system provides a natural bridge between the reactivity of the original system (in the limit of $p=1$) and Lyapunov exponents of the original system (in the limit of $p\to \infty$), in the following section, we present an example that illustrates the implications and advantages of studying the reactivity of the $p$-iteration system.

\section{Application to Network Synchronization}\label{sec:sync}

Consider a network of $n$ coupled 
maps where the state of map $i$ 
evolves over time according to the following equation:
\begin{equation} \label{eq:synch}
    \bx_{k+1}^i = \bF(\bx_k^i) - \kappa \sum_j L_{ij} \bH (\bx_k^j), \quad i = 1, \hdots, n.
\end{equation}
Here, 
$\bx_k^i \in \mathbb{R}^m$ is the state vector \ZA{of map $i$ at time $k$}, $\bF: \mathbb{R}^m \rightarrow \mathbb{R}^m$ describes the dynamics of the map and $\bH: \mathbb{R}^m \rightarrow \mathbb{R}^m$ is the coupling function.
The scalar $\kappa \geq 0$ is the coupling strength. 
The network connectivity is described by the adjacency matrix $A=[A_{ij}]$. For simplicity, we assume the graph is undirected, so the matrix $A$ is symmetric. The entry $A_{ij} = A_{ji} \geq 0$ is the strength of the coupling between nodes $i$ and $j$.
The corresponding Laplacian matrix $L = [L_{ij}]$ is defined as $L_{ij}=-A_{ij}, \ i \neq j$, and $L_{ii}=\sum_{j} A_{ij}$.
Since $L$ is symmetric, it is diagonalizable via an orthogonal transformation. Let $T\in \mathbb{R}^{n\times n}$ be an orthogonal matrix whose columns are the eigenvectors of $L$. Then, 
$T^\top L T = \Lambda$ where $\Lambda = \text{diag}(0,\lambda_2, \lambda_3, \hdots, \lambda_n) \in \mathbb{R}^{n\times n}$ and $0 < \lambda_2 \leq \lambda_3 \leq \cdots \leq \lambda_n$ are the eigenvalues of $L$. The first column of $T$ is $\hat{\pmb{1}}=1/\sqrt{n} [1, 1, \hdots,1]^\top$ which corresponds to the zero eigenvalue of $L$. 

The synchronous solution \ZA{of \eqref{eq:synch} is defined by $\{\bx_k^1=\ldots=\bx_k^n=\bs_k\}$ and governed by}
   $ \bs_{k+1} = \bF (\bs_k).$
An attractor $\mathcal{A}$ is a set in $\mathbb{R}^m$ toward which the sequence $ \bs_k$ converges as $k \to \infty$.
An attractor can be either a fixed point, a limit cycle, or a chaotic attractor.
There can be more than one attractor. Our analysis is attractor-specific, i.e., it depends on the specific attractor on which the dynamics converges.
The initial conditions are chosen randomly on the attractor $\mathcal{A}$. Then $\bs_{k+1} = \bF (\bs_k)$ is evaluated for a long enough time until the orbit converges on the attractor.

Next, we study the linear stability of the synchronous solution.
We define the infinitesimal perturbations about the synchronous solution as $\delta \bx^i_k := \bx_k^i - \bs_k$.
The time evolution of these perturbations is described by
$
    \delta \bx^i_{k+1} = D \bF(\bs_k) \delta \bx^i_k - \kappa \sum_j L_{ij} D \bH (\bs_k) \delta \bx^j_k,
$
where $D \bF$ ($D\bH$) is the Jacobian of $\bF$ ($\bH$). 
This equation can be written in compact form,
\begin{equation} \label{eq:compact0}
    \delta\bX_{k+1} = \left( I_n \otimes D\bF(\bs_k) - \kappa L \otimes D \bH(\bs_k) \right) \delta \bX_k,
\end{equation}
where $\delta \bX_k = \left[{\delta \bx^1_k}^\top,\hdots,  {\delta \bx^n_k}^\top\right]^\top$ and $\otimes$ is the Kronecker product.

We define the matrix $ V \in \mathbb{R}^{n \times (n-1)}$ as the matrix formed by all columns of $T$ except the first, i.e., $T = [\hat{\pmb{1}}, V]$.
By pre-multiplying \eqref{eq:compact0} by $V^\top \otimes I_m$, see \cite{pecora1998master}, the equation for the time evolution of the perturbations transverse to the synchronization manifold is obtained:
\begin{equation} \label{eq:compact}
    \delta \hat{\bX}_{k+1} = Z(\bs_k) \delta \hat{\bX}_k,
\end{equation}
where $Z(\bs_k) = I_{n-1} \otimes D\bF(\bs_k) - \kappa \hat{\Lambda} \otimes D \bH(\bs_k)$.
Here, $\hat{\Lambda} = \text{diag}(\lambda_2, \lambda_3, \hdots, \lambda_n) = V^\top L V \in \mathbb{R}^{n-1 \times n -1}$.
Then, the local stability of the synchronization manifold is determined through the stability of \eqref{eq:compact}.

We consider an example for which we take the chaotic Henon map, with $\bx = [x, \ y]^\top$ and the 
function and the Jacobian
\begin{equation} \label{eq:fh}
    \bF(\bx) = \begin{bmatrix}
        y + 1 - 1.4 x^2 \\
        0.3 x
    \end{bmatrix},  D\bF (\bx) = \begin{bmatrix}
        -2.8 x & 1 \\
        0.3 & 0
    \end{bmatrix},
\end{equation}
and $\bH(\bx) = \bF(\bx)$.
This indicates that
\begin{equation} \label{eq:Zhenon}
    Z(\bs_k) = \left(I_{n-1} - \kappa \hat{\Lambda} \right) \otimes D \bF(\bs_k).
\end{equation}
We now study the reactivity of the $p$-iteration system of \eqref{eq:compact}:
$\delta \bX_{k+p} = Z(\bs_k)^{(p)} \delta \bX_k$,
where
\begin{align} \label{eq:piterz}
\begin{split}
    Z(\bs_k)^{(p)} & = Z(\bs_{k+p-1}) \cdots Z(\bs_{k+1}) Z(\bs_{k}) \\
    & = \left(I_{n-1} - \kappa \hat{\Lambda} \right)^p \otimes  D\bF(\bs_k)^{(p)},
\end{split}
\end{align}
where 
$D\bF(\bs_k)^{(p)} = D\bF(\bs_{k+p-1}) 
\allowbreak \cdots 
\allowbreak D\bF(\bs_{k})$. 
Note that $Z(\bs_k)^{(1)} = Z(\bs_k)$.

In what follows, we let the norm $\|\cdot\|$ be a $Q-$weighted $L_{\mathfrak{p}}$ 
where $Q$ is a positive diagonal
matrix. 
The reactivity of the $p$-iteration linearized map is 
\begin{equation*}
    r_{\| \cdot \|}[Z(\bs_k)^{(p)}]=\left \| \left(I_{n-1} - \kappa \hat{\Lambda} \right)^p \right\|  \left\| D\bF(\bs_k)^{(p)} \right\| -1.
\end{equation*}
The term $ \| (I_{n-1} - \kappa \hat{\Lambda} )^p \|$ describes the contribution of the topology to reactivity, and the term $\left\| D\bF(\bs_k)^{(p)}  \right\|$ denotes the effect of dynamics. 
The sup-mean and the lim-mean reactivities of \eqref{eq:compact} are
$ \limR (p) = \lim_{N \to \infty} 1/N \sum_{k=0}^{N-1} r_{\| \cdot \|}[Z(\bs_k)^{(p)} ]$, and $\barR (p) = \allowbreak \sup_{N} \allowbreak 1/N \sum_{k=0}^{N-1} r_{\| \cdot \|}[Z(\bs_k)^{(p)} ]$,
respectively.
Since points on the chaotic synchronous solution $\bs_k \in \mathcal{A}, \forall k$, $\limR (p)$ and $\barR (p)$ become
\begin{subequations}
\begin{align}
    \limR {(p)} & = \left \langle r_{\| \cdot \|}[Z(\bs)^{(p)} ] \right \rangle_{\bs \in \mathcal{A}}, \\
    \barR {(p)} & = \max_{\bs \in \mathcal{A}} r_{\| \cdot \|}[Z(\bs)^{(p)} ], 
\end{align}
\end{subequations}
where $\langle \cdot \rangle_{\bs \in \mathcal{A}}$ denotes average over all points on the attractor $\mathcal{A}$.
We define positive attractor-specfic constants $\beta_\mean^{(p)} := \left \langle \left\| D\bF(\bs)^{(p)} \right\| \right \rangle_{\bs \in \mathcal{A}}$ and $\beta_{\max}^{(p)} := \max_{\bs \in \mathcal{A}} \left\| D\bF(\bs)^{(p)} \right\|$.
It trivially follows $\beta_{\max}^{(p)} \geq \beta_\mean^{(p)}$.
As a result, the lim-mean and the sup-mean reactivities become $\limR {(p)} = \beta_\mean^{(p)} \| (I_{n-1} - \kappa \hat{\Lambda})^p \| -1$ and $\barR {(p)} = \beta_{\max}^{(p)} \| (I_{n-1} - \kappa \hat{\Lambda} )^p\| -1$.

Next, we find the range of the coupling strength $\kappa$ in which the transverse dynamics to the synchronization manifold are contractive. 
The range of $\kappa_{\min} < \kappa < \kappa_{\max}$ calculated based on the sup-mean reactivity $\barR {(p)}$ is expected to be a subset of the range calculated via the lim-mean reactivity $\limR {(p)}$ since $\limR {(p)}$ requires more relaxed conditions (since $\beta_{\max}^{(p)} \geq \beta_\mean^{(p)}$). 
Recalling that $I_{n-1} - \kappa \hat{\Lambda} = \text{diag} (1-\kappa\lambda_2, 1-\kappa\lambda_3, \hdots, 1-\kappa \lambda_n)$, it follows $\| (I_{n-1} - \kappa \hat{\Lambda})^p \| = \max_{i=2,3,\hdots,n} | (1 - \kappa \lambda_i)^p|$. 
We enforce $\limR {(p)} < 0$ which results in $\max_{i=2,3,\hdots,n} | (1 - \kappa \lambda_i)^p| < 1/{\beta_\mean^{(p)}}$. 
It follows $(1-\kappa_{\min} \lambda_2)^p < 1/{\beta_\mean^{(p)}}$ and $(\kappa_{\max} \lambda_n -1 )^p< 1/{\beta_\mean^{(p)}}$. After simplification
\begin{equation} \label{eq:kappa}
    \mathcal{L}_\mean^{(p)}:=\dfrac{{\beta_\mean^{(p)}}^{1/p} - 1}{\lambda_2 {\beta_\mean^{(p)}}^{1/p}} < \kappa < \dfrac{{\beta_\mean^{(p)}}^{1/p}+1}{\lambda_n {\beta_\mean^{(p)}}^{1/p}}=:\mathcal{U}_\mean^{(p)},
\end{equation}
where $\mathcal{L}_\mean^{(p)}$ and $\mathcal{U}_\mean^{(p)}$ are lower and upper bounds.
Note that $\beta_\mean^{(p)} > 1$ for chaotic systems (here the chaotic Henon map) since these systems have a positive Maximum Lyapunov Exponent.
The network is lim-mean synchronizable, i.e., $\exists \kappa$ such that $\limR (p) <0 $ if
\begin{equation} \label{eq:synchronizabilitymean}
    \dfrac{\lambda_n}{\lambda_2} < \dfrac{{\beta_\mean^{(p)}}^{1/p}+1}{{\beta_\mean^{(p)}}^{1/p}-1} := \mathcal{S}_{\mean}^{(p)},
\end{equation}
where $\mathcal{S}_{\text{mean}}^{(p)}$ is the synchronizability based on $\limR$.
Similar derivations to \eqref{eq:kappa}-\eqref{eq:synchronizabilitymean} can be obtained by enforcing $\barR (p)< 0$ instead of $\limR < 0$, where in that case, $\beta_\mean^{(p)}$ must be replaced by $\beta_{\max}^{(p)}$.
It follows
\begin{equation} \label{eq:synchronizabilitymax}
    \dfrac{\lambda_n}{\lambda_2} < \dfrac{{\beta_{\max}^{(p)}}^{1/p}+1}{{\beta_{\max}^{(p)}}^{1/p}-1} := \mathcal{S}_{\max}^{(p)},
\end{equation}
We also denote $\mathcal{L}_{\max}^{(p)}$ and $\mathcal{U}_{\max}^{(p)}$ as the lower and the upper bounds of $\kappa$, respectively, when $\beta_{\max}^{(p)}$ is used in \eqref{eq:kappa} instead of $\beta_\mean^{(p)}$.
First, the bounds for the coupling strength $\kappa$ are tighter with $\barR (p) < 0$ since for a fixed spectrum, $\beta_{\max}^{(p)} \geq \beta_\mean^{(p)}$, so the upper bound $\mathcal{U}_{\max}^{(p)} \leq \mathcal{U}_\mean^{(p)}$. 
This is easily verifiable since the function $g(y) = (y+1)/y$ is strictly decreasing for $y > 0$. 
Similarly, the lower bound $\mathcal{L}_\mean^{(p)} \leq \mathcal{L}_{\max}^{(p)}$.
This is also easily verifiable since the function $h(y) = (y-1)/y$ is strictly increasing for $y>0$.
For sup-mean synchronizability, i.e., $\exists \kappa$ such that $\limR (p) <0 $, it is easy to verify that $\mathcal{S}_{\max}^{(p)} \leq \mathcal{S}_{\mean}^{(p)}$.
Therefore, all conditions on the sup-mean approach are tighter than those of the lim-mean approach.

\begin{figure}
    \centering
    \includegraphics[width=\linewidth]{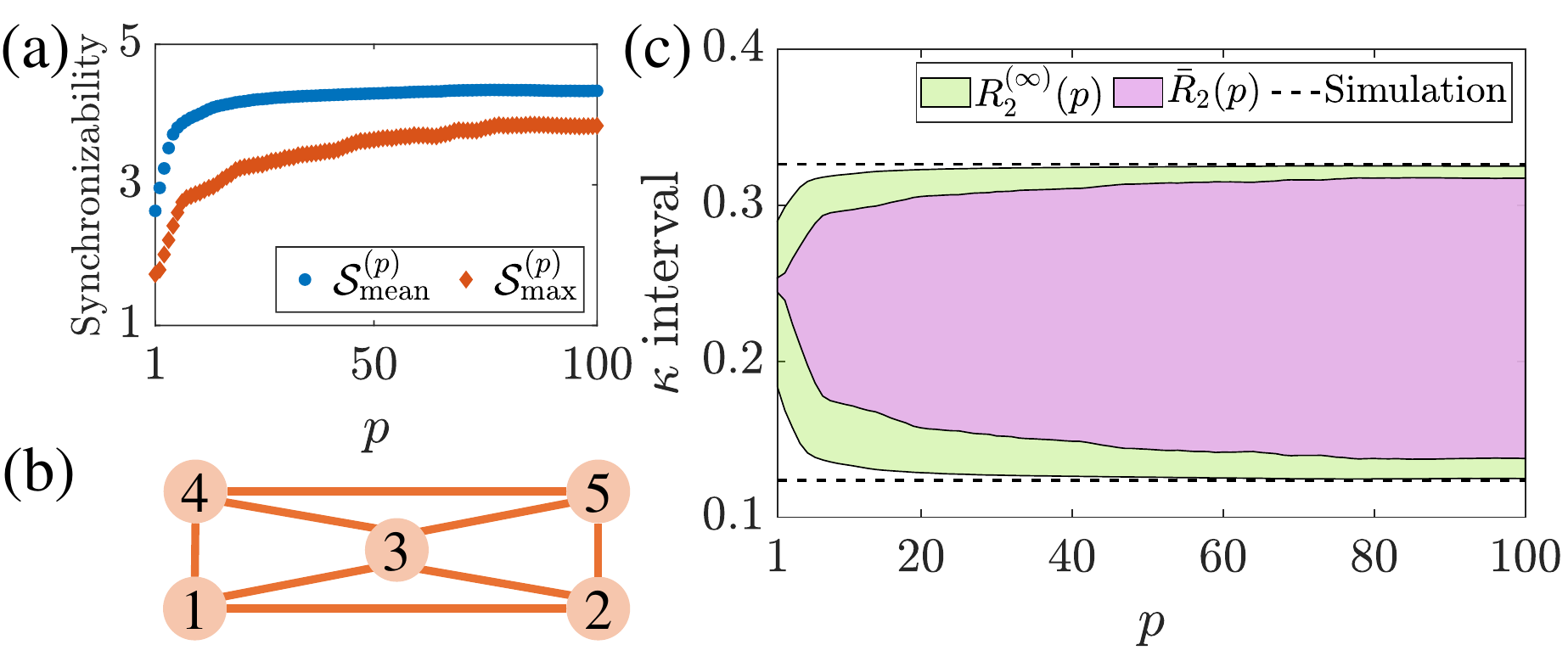}
    \caption{ (a) Synchronizability criteria based on \eqref{eq:synchronizabilitymean}, computed using either $R_2^{(\infty)} (p)$ $\left(\beta_\mean^{(p)}\right)$ or $\bar{R}_2 (p)$ $\left(\beta_{\max}^{(p)}\right)$.  
    (b) The unweighted and undirected network used in this example.  
    (c) Bounds of the synchronous coupling strength interval $\kappa$, as defined in \eqref{eq:kappa}. The green region corresponds to the interval based on $R^{(\infty)}_2(p)$, while the pink region is based on $\bar{R}_2(p)$.  
    Dashed black lines indicate the interval of $\kappa$ in which the synchronization of the full system \eqref{eq:synch} is stable. Synchronization is considered achieved when after simulating \eqref{eq:synch}, the error
  $ 
    E = \frac{1}{k_f - k_0} \sum_{k=k_0}^{k_f} \frac{1}{n} \sum_{i=1}^{n} \big\| \bx^i_k - \frac{1}{n} \sum_{j=1}^{n} \bx^j_k \big\|_2
  $  falls below $10^{-10}$. In this example, we set $k_f = 10000$ and $k_0 = 9000$.}
\label{fig:henon}
\end{figure}

Next, we show the effect of $p$ on the bounds predicted by $\bar{R}_2 (p)$ and $R^{(\infty)}_2 (p)$ via a numerical example when $L_2$ norm is used.
We first plot the synchronizability measures for the chaotic Henon map in \eqref{eq:synchronizabilitymean} and \eqref{eq:synchronizabilitymax} in Fig.~\ref{fig:henon}(a). 
Note that for all values of $p\in[1 \;100]$, the synchronizability $\mathcal{S}_{\mean}^{(p)}$ is higher than that based on $\mathcal{S}_{\max}^{(p)}$.
This suggests that a network of coupled Henon maps may be synchronizable based on $R^{(\infty)}_2 (p)$, but not synchronizable based on $\bar{R}_2 (p)$.
Also, as the number of iterations $p$ grows, the synchronizabilities grow.

We take a 5-node unweighted and undirected network shown in Fig.\,\ref{fig:henon}(b).
The Laplacian matrix $L$ of this network has eigenvalues $\lambda_1 = 0$, $\lambda_2 = \lambda_3 = 3$, and $\lambda_4 = \lambda_5 = 5$.
Note that $\lambda_5/\lambda_2 \approx 1.67$, which is lower than the minimum synchronizability required in Fig.\,\ref{fig:henon}(a); thus, the network is synchronizable for all $p$-iterations.
In Fig.\,\ref{fig:henon}(c), we plot the bounds from \eqref{eq:kappa} based on the choice of both $R^{(\infty)}_2 (p)$ and $\bar{R}_2 (p)$ measures.
As expected, the interval based on $\bar{R}_2 (p)$ contains inside interval based on $R^{(\infty)}_2 (p)$. 
Also, as $p$ grows, the upper and lower bounds of the intervals increase and decrease, respectively, until for large $p$, they converge to the synchronous bounds predicted by numerically simulating \eqref{eq:synch}.

\color{black}

\section{Conclusion and future work}
\label{sec:conc}
The closely related concepts of reactivity and contraction have been independently introduced and investigated in the literature. Here we present a unified framework to study reactivity and contractility, and connect them with the concept of finite-time Lyapunov exponents. 
We focus on discrete-time systems and their $p$-iterations where $p>1$. 
    We introduce a definition of reactivity for nonlinear discrete-time systems based on general 
    norms. 
    Stability conditions based on reactivity are provided for both time-varying and time-invariant nonlinear systems. The case of linear time-varying systems is relevant to the study of the stability of a particular trajectory under the assumption of small perturbations, which we use to characterize the local stability of the synchronous solution for networks of coupled maps. For this case, we can replace the traditional requirement for stability that the trajectory is never reactive with a less restrictive requirement that the trajectory is non-reactive `on average'. Also, for this case, the reactivity of the $p$-iteration system coincides with the finite-time (time equal $p$) Lyapunov exponent. In the limit of infinite $p$, the concepts of reactivity and Lyapunov exponents coincide, which allows us to bridge the gap between contraction theory and the master stability function theory used to assess the stability of the synchronous solution in networks.

    A topic that is left for future investigation is the determination of a stopping criterion for increasing values of the iteration number $p$ to use in our propositions.

We emphasize that our work introduces a tradeoff between computational expenditure and the advantages of the analysis.  In particular, we see that as we increase the number of iterations $p$, computation of the contraction region becomes more expensive, but we are able to predict a wider range of stability. This can be seen in Fig.~\ref{fig:henon} for the case of the stability of the synchronous solution for a network of coupled nonlinear maps.


\begin{thebibliography}{10}

\bibitem{abarbanel1991lyapunov}
Henry~DI Abarbanel, Reggie Brown, and MB~Kennel.
\newblock Lyapunov exponents in chaotic systems: their importance and their
  evaluation using observed data.
\newblock {\em International Journal of Modern Physics B}, 5(09):1347--1375,
  1991.

\bibitem{aminzare2013logarithmic}
Zahra Aminzare and Eduardo~D Sontag.
\newblock Logarithmic lipschitz norms and diffusion-induced instability.
\newblock {\em Nonlinear Analysis: Theory, Methods \& Applications}, 83:31--49,
  2013.

\bibitem{aminzare2014contraction}
Zahra Aminzare and Eduardo~D Sontag.
\newblock Contraction methods for nonlinear systems: A brief introduction and
  some open problems.
\newblock In {\em 53rd IEEE Conference on Decision and Control}, pages
  3835--3847. IEEE, 2014.

\bibitem{Asllani2018Structure}
Malbor Asllani, Renaud Lambiotte, and Timoteo Carletti.
\newblock Structure and dynamical behavior of non-normal networks.
\newblock {\em Science Advances}, 4(12):eaau9403, 2018.

\bibitem{bessaga1959converse}
Cz~Bessaga.
\newblock On the converse of the banach fixed point principle.
\newblock In {\em Colloq. Math.}, volume~7, pages 41--43, 1959.

\bibitem{Biancalani2017Giant}
Tommaso Biancalani, Farshid Jafarpour, and Nigel Goldenfeld.
\newblock Giant amplification of noise in fluctuation-induced pattern
  formation.
\newblock {\em Phys. Rev. Lett.}, 118:018101, Jan 2017.

\bibitem{bullo2022contraction}
Francesco Bullo.
\newblock {\em Contraction theory for dynamical systems}.
\newblock Kindle Direct Publishing, 1.1 edition, 2022.

\bibitem{Duan2022Network}
Chao Duan, Takashi Nishikawa, Deniz Eroglu, and Adilson~E. Motter.
\newblock Network structural origin of instabilities in large complex systems.
\newblock {\em Science Advances}, 8(28):1--12, 2022.

\bibitem{Farrell1996Generalized}
Brian~F. Farrell and Petros~J. Ioannou.
\newblock Generalized stability theory. part i: Autonomous operators.
\newblock {\em Journal of Atmospheric Sciences}, 53(14):2025 -- 2040, 1996.

\bibitem{Granas2003}
Andrzej Granas and James Dugundji.
\newblock {\em Elementary Fixed Point Theorems}, pages 9--84.
\newblock Springer New York, New York, NY, 2003.

\bibitem{Gudowska2020From}
E.~Gudowska-Nowak, M.~A. Nowak, D.~R. Chialvo, J.~K. Ochab, and W.~Tarnowski.
\newblock {From Synaptic Interactions to Collective Dynamics in Random Neuronal
  Networks Models: Critical Role of Eigenvectors and Transient Behavior}.
\newblock {\em Neural Computation}, 32(2):395--423, 02 2020.

\bibitem{Hennequin2012Nonnormal}
Guillaume Hennequin, Tim~P. Vogels, and Wulfram Gerstner.
\newblock Non-normal amplification in random balanced neuronal networks.
\newblock {\em Phys. Rev. E}, 86:011909, Jul 2012.

\bibitem{Lindmark2021Centrality}
Gustav Lindmark and Claudio Altafini.
\newblock Centrality measures and the role of non-normality for network control
  energy reduction.
\newblock {\em IEEE Control Systems Letters}, 5(3):1013--1018, 2021.

\bibitem{lohmiller1998contraction}
Winfried Lohmiller and Jean-Jacques~E Slotine.
\newblock On contraction analysis for non-linear systems.
\newblock {\em Automatica}, 34(6):683--696, 1998.

\bibitem{lohmiller2000nonlinear}
Winfried Lohmiller and Jean-Jacques~E Slotine.
\newblock Nonlinear process control using contraction theory.
\newblock {\em AIChE journal}, 46(3):588--596, 2000.

\bibitem{MUOLO2019Patterns}
Riccardo Muolo, Malbor Asllani, Duccio Fanelli, Philip~K. Maini, and Timoteo
  Carletti.
\newblock Patterns of non-normality in networked systems.
\newblock {\em Journal of Theoretical Biology}, 480:81--91, 2019.

\bibitem{nazerian2024efficiency}
Amirhossein Nazerian, Joseph~D Hart, Matteo Lodi, and Francesco Sorrentino.
\newblock The efficiency of synchronization dynamics and the role of network
  syncreactivity.
\newblock {\em Nature Communications}, 15(1):9003, 2024.

\bibitem{nazerian2022matryoshka}
Amirhossein Nazerian, Shirin Panahi, Ian Leifer, David Phillips, Hernan Makse,
  and Francesco Sorrentino.
\newblock Matryoshka and disjoint cluster synchronization of networks, 2022.

\bibitem{nazerian2023commphys}
Amirhossein Nazerian, Shirin Panahi, and Francesco Sorrentino.
\newblock Synchronization in networked systems with large parameter
  heterogeneity.
\newblock {\em Communications Physics}, 6(1):253, 2023.

\bibitem{Nazerian2023epl}
Amirhossein Nazerian, Shirin Panahi, and Francesco Sorrentino.
\newblock Synchronization in networks of coupled oscillators with mismatches.
\newblock {\em Europhysics Letters}, 143(1):11001, jul 2023.

\bibitem{nazerian2023reactability}
Amirhossein Nazerian, David Phillips, Mattia Frasca, and Francesco Sorrentino.
\newblock The reactability of discrete time systems.
\newblock {\em IEEE Control Systems Letters}, 7:3657--3662, 2023.

\bibitem{nazerian2023single}
Amirhossein Nazerian, David Phillips, Hern{\'a}n~A Makse, and Francesco
  Sorrentino.
\newblock Single-integrator consensus dynamics over minimally reactive
  networks.
\newblock {\em IEEE Control Systems Letters}, 2023.

\bibitem{Neubert1997ALTERNATIVES}
Michael~G. Neubert and Hal Caswell.
\newblock Alternatives to resilience for measuring the responses of ecological
  systems to perturbations.
\newblock {\em Ecology}, 78(3):653--665, 1997.

\bibitem{ott2002chaos}
Edward Ott.
\newblock {\em Chaos in dynamical systems}.
\newblock Cambridge university press, 2002.

\bibitem{panahi2021cluster}
Shirin Panahi, Isaac Klickstein, and Francesco Sorrentino.
\newblock Cluster synchronization of networks via a canonical transformation
  for simultaneous block diagonalization of matrices.
\newblock {\em Chaos: An Interdisciplinary Journal of Nonlinear Science},
  31(11):111102, 2021.

\bibitem{panahi2022pinning}
Shirin Panahi, Matteo Lodi, Marco Storace, and Francesco Sorrentino.
\newblock Pinning control of networks: Dimensionality reduction through
  simultaneous block-diagonalization of matrices.
\newblock {\em Chaos: An Interdisciplinary Journal of Nonlinear Science},
  32(11), 2022.

\bibitem{panahi2021group}
Shirin Panahi and Francesco Sorrentino.
\newblock Group synchrony, parameter mismatches, and intragroup connections.
\newblock {\em Physical Review E}, 104(5):054314, 2021.

\bibitem{pecora1998master}
Louis~M Pecora and Thomas~L Carroll.
\newblock {Master stability functions for synchronized coupled systems}.
\newblock {\em Physical review letters}, 80(10):2109, 1998.

\bibitem{slotine2001modularity}
J-JE Slotine and Winfried Lohmiller.
\newblock Modularity, evolution, and the binding problem: a view from stability
  theory.
\newblock {\em Neural networks}, 14(2):137--145, 2001.

\bibitem{slotine2003modular}
Jean-Jacques~E Slotine.
\newblock Modular stability tools for distributed computation and control.
\newblock {\em International Journal of Adaptive Control and Signal
  Processing}, 17(6):397--416, 2003.

\bibitem{Tang2014Reactivity}
Si~Tang and Stefano Allesina.
\newblock Reactivity and stability of large ecosystems.
\newblock {\em Frontiers in Ecology and Evolution}, 2, 2014.

\bibitem{wang2005partial}
Wei Wang and Jean-Jacques~E Slotine.
\newblock {On partial contraction analysis for coupled nonlinear oscillators}.
\newblock {\em Biological Cybernetics}, 92(1):38--53, jan 2005.

\bibitem{Wolf1986Quantifying}
A.~Wolf.
\newblock {\em 13. Quantifying chaos with Lyapunov exponents}, pages 273--290.
\newblock Princeton University Press, Princeton, 1986.

\bibitem{Zhou2017Asymptotic}
Bin Zhou and Tianrui Zhao.
\newblock On asymptotic stability of discrete-time linear time-varying systems.
\newblock {\em IEEE Transactions on Automatic Control}, 62(8):4274--4281, 2017.

\end{thebibliography}

\end{document}